\documentclass[a4paper,12pt,twoside]{article}
\usepackage{inputenc}
\usepackage[english,russian]{babel} 
\usepackage{jcem}
\usepackage{amsfonts}
\usepackage{amsmath}
\usepackage{amssymb}
\usepackage{amsthm}
\usepackage{graphicx}
\newcommand{\mydef}{\operatorname{:\!=}}
\newcommand{\fix}[1]{\ensuremath{{\mathbf{Fix}({#1})}}}
\newcommand{\myle}{\ensuremath{\mathrel\preccurlyeq}}                
\newcommand{\myles}{\mathrel\prec}                          
\newcommand{\myor}{\ensuremath{\lor}}        
\newcommand{\myemp}{\ensuremath{\varnothing}}   
\newcommand{\myimp}{\ensuremath{\mathrel\Rightarrow}}
\newcommand{\myll}{\ensuremath{\forall}}         
\newcommand{\beq}{\begin{equation}}
\newcommand{\beqnt}{\begin{equation}\notag}
\newcommand{\eeq}{\end{equation}}
\newcommand{\bml}{\begin{multline}}
\newcommand{\eml}{\end{multline}}
\newcommand{\bmlnt}{\begin{multline*}}
\newcommand{\emlnt}{\end{multline*}}

\newcommand{\bpra}[1]{\begin{proof}[#1]}

\newcommand{\fref}[1]{{\rm(\ref{#1})}}                      
\newcommand{\NA}{{\mathbb N}}

\newcommand{\iter}[1]{\ensuremath{\mathrm{Iter({#1})}}}
\newcommand{\ext}[1]{\ensuremath{\mathrm{Ext({#1})}}}
\newcommand{\ffam}{\ensuremath{\mathcal{F}}}

\begin{document}
\def\Re{\mathop{\rm Re}\,}
\def\B{\mathop{\frak B}}
\def\R{\mathop{\Bbb R}}
\def\A{\mathop{\frak A}}
\def\G{\mathop{\bf G}}
\def\U{\mathop{\frak U}}
\def\H{\mathop{\cal H}}
\def\Im{\mathop{\rm Im}\,}
\def\dom{\mathop{\rm dom}\,}
\def\dist{\mathop{\rm dist}}
\def\grad{\mathop{\rm grad}}
\renewcommand{\proof}{\vspace{2mm}\hspace{-7mm}\textit{Proof.}}
\renewcommand{\endproof}{\begin{flushright} \vspace{-2mm}$\Box$\vspace{-4mm}
\end{flushright}}
\newcommand{\phan}{\hspace*{0cm}}
\newcommand{\comment}{}

\begin{flushleft}
\textbf{MSC 47H10, 54C10, 54E45, 91B50} 
\end{flushleft}
\author{Serkov Dmitrii Aleksandrovich, {\rm Krasovskii Institute of Mathematics and Mechanics, 
Ural Federal University named after B.N.Yeltsin,
Yekaterinburg, Russian Federation, e-mail: serkov@imm.uran.ru}}
\title{On approximation of joint fixed points}
\maketitle{}

\begin{abstract} \begin{tabular}{p{0mm}p{139mm}}
&\noindent {\footnotesize \qquad 
For a given partially ordered set (poset) and a given family of mappings of the poset into itself, we study the problem of the description of joint fixed points of this family. Well-known Tarski's theorem gives the structure of the set of joint fixed points of isotone automorphisms on a complete lattice. This theorem has several generalizations (see., e.g., Markowsky, Ronse) that weaken demands on the order structure and upgrade in an appropriate manner the assertion on the structure of the set of joint fixed points. However, there is a lack of the statements similar to Kantorovich or Kleene theorems, describing the set of joint fixed points in terms of convergent sequences of the operator degrees. The paper provids conditions on the poset and on the family; these conditions ensure that the iterative sequences of elements of this family approximate the set of joint fixed points. The result obtained develops in a constructive direction the mentioned theorems on joint fixed points.

\qquad}
\end{tabular}\end{abstract}

\markboth{}{}

\section*{Introduction}
\hspace{0.7 cm}
The key object in the solution of the differential game of quality --- the maximal stable bridge \cite{KraSub88e} --- is the greatest fixed point of a special operator (programmed absorption operator) acting in the boolean of the extended phase space of controlled system \cite{Chentsov76MS}. In its turn, the action of this operator is described by actions of some operator family, in such a way that any joint fixed point of this family is a fixed point of the original operator. This relation makes it possible to describe the object of interest --- the maximal stable bridge --- in terms of the operator family with a relatively simple structure.

The paper presents the theorem on the joint fixed points of a family of mappings arising in connection with the above. This result develops in a constructive direction well-known Tarski's theorem \cite[theorem 2]{Tarski1955} on the structure of the set of joint fixed points of automorphisms on a complete lattice. This theorem has a number of generalizations (see, e.g., \cite[theorem 10]{Markow1976}, \cite[corollary 3.3]{Ronse1994}),  weakening demands on the order structure and upgrading in an appropriate manner the assertion on the structure of the set of joint fixed points. However, there is a lack of known statements like the Kantorovich theorem (see. \cite{Kantorovitch1939AM}) or Kleene theorem that describe the set of joint fixed points in terms of convergent sequences of the operator degrees. theorem 1 below fills this gap. 

\section{Results}
\hspace{0.7 cm}

Let $(X, \myle)$ be a partially ordered set (poset). A set $C\subset X$ is called a \emph{chain} if it is totally ordered by \myle: 
$$
(x\myle y)\myor(y\myle x)\qquad \myll x,y\in C.
$$ 
In particular, \myemp\ is a chain. Following \cite{Markow1976} we call a poset $(X,\myle)$ \emph{chain--complete} poset if, for every chain $C\subset X$, there exists the least upper bound $\sup C$ of the chain $C$ in $X$. In particular, a chain--complete poset $(X, \myle)$ contains the least element $\bot\in X$ (as the least upper bound of the empty chain), and thus it is not empty. Let $(W, \myles)$ and $(X, \myle)$ be chain complete posets and $f:W \mapsto X$. The mapping $f$ is called \emph{isotone} if $\myll x,y\in X$ $(x\myles y)\myimp(f(x)\myle f(y))$. We say that the mapping $f$ is \emph{chain--continuous} if, for any non--empty chain $C\subset W$, the equality $f(\sup C)=\sup\{f (w): w\in C\}$ holds. We say that $(X, \myle)$ is a \emph{complete} poset \cite[p.1.2.1]{BAREN1981E}, if there is the least element $\bot\mydef\sup\myemp$ and each directed set $D\subset X$ has the least upper bound $\sup D$. Recall that a set is called \emph{directed} if any its finite subset has a majorant.

Let $ X\neq\myemp$ and a non--empty family \ffam\ of mappings $f:X\mapsto X$ $\myll f\in\ffam$ be given. We call the set $\ext{\ffam}\mydef\{x\in X\mid x\myle f(x)\ \myll f\in\ffam\}$ the \emph{extensivity domain} of the family \ffam. Denote the set of all finite compositions of mappings from \ffam\ by
$$
\iter{\ffam}\mydef\{fg\ldots h\mid f,g,\ldots,h\in\ffam\};
$$
here $fg\ldots h$ stands for the composition of mappings $f$,$g$,\ldots, $h$:
$$
fg\ldots h(x)\mydef f(g(\ldots h(x)\ldots))\qquad\myll x\in X.
$$
If $f=g=\ldots=h$, then we denote $f^k\mydef fg\ldots h$, where $k\mydef |\{f,g,\ldots, h\}|$ is the number of elements in the composition $fg\ldots h$.
Thus, $\iter{\ffam}\subset{X^X}$. We say that the family \ffam\ is \emph{commutative} if
$$
fg=gf\qquad\myll f,g \in\ffam.
$$
Denote by $\fix{f}\subset X$ the \emph{set of fixed points} of mapping $f\in\ffam$: $\fix{f}\mydef\{x \in X\mid x=f(x)\}$. Finally, we denote by $\fix{\ffam}\subset X$ the \emph{set of joint fixed points} of the family \ffam: $\fix{\ffam}\mydef\bigcap_{f \in\ffam}\fix{f}$, i.e. $f(x)=x$ $\myll x\in\fix{\ffam}$ $\myll f\in\ffam$.

\begin{theorem}\label{teo-limit-eq-fix-3} 
Let $(X,\myle)$ be a chain--complete poset and let \ffam\ be a commutative family of isotone mappings from $X$ into itself. Then 

\textup{(i)} if $X$ is a complete lattice, then $\fix{\ffam}$ is a complete lattice with respect to the induced order; thus it has the least element $\bot_{\fix{\ffam}}$ and, in particular, $\fix{\ffam}\neq\myemp$ \textup{(}Tarski\textup{)}; 

\textup{(ii)} $\fix{\ffam}$ is chain--complete with respect to the induced order; thus it has the least element $\bot_{\fix{\ffam}}$ and, in particular, $\fix{\ffam}\neq\myemp$ \textup{(}Markowsky\textup{)};

\textup{(iii)} if elements of \ffam\ are chain--continuous, then
\beq\label{fixfam}
\fix{\ffam}=\{\sup\{\varphi (x)\mid\varphi\in\iter{\ffam}\}\mid x\in\ext{\ffam}\},
\eeq 
and, in particular, 
\beq\label{botfixfam}
\bot_{\fix{\ffam}}=\sup\{\varphi(\bot)\mid\varphi\in\iter{\ffam}\}.
\eeq 
\end{theorem}

In the case when the family \ffam\ is a singleton, we obtain the corollary.

\begin{corollary}\label{teo-limit-eq-fix-2}
Let $f:X\mapsto X$  be an isotone mapping on a chain--complete poset $(X,\myle)$. Then 

\textup{(i)} if $X$ is a complete lattice, then $\fix{f}$ is a complete lattice with respect to the induced order; thus so it has the least element $\bot_{\fix{f}}$ and, in particular, $\fix{f}\neq\myemp$ \textup{(}Tarski\textup{)}, 

\textup{(ii)} $\fix{f}$ is chain--complete with respect to the induced order; thus it has the least element $\bot_{\fix{\ffam}}$ and, in particular, $\fix{f}\neq\myemp$ \textup{(}Markowsky\textup{)},

\textup{(iii)} if $f$ is chain--continuous, then \textup{(}Kleene\textup{)}
\beq\label{botfixf}
\bot_{\fix{f}}=\sup\{f^n(\bot)\mid n\in\NA\},
\eeq 

\textup{(iv)} if $f$ is chain--continuous, then 
\beq\label{limit-eq-fix-2}
\fix{f}=\{\sup\{f^n(x)\mid n\in\NA\}\mid x\myle f(x)\}.
\eeq
\end{corollary}


\section{Proofs}

Assertion (i) of theorem \ref{teo-limit-eq-fix-3} is contained in \cite[theorem 2]{Tarski1955}, assertion (ii) of theorem \ref{teo-limit-eq-fix-3} is contained in \cite[theorem 10]{Markow1976}. Let us turn to the proof of the assertion (iii) of theorem \ref{teo-limit-eq-fix-3}.

First, we note that, by virtue of \cite[p.33]{Cohn1965} (see also \cite[corollary 2]{Markow1976}) and of chain--completeness condition, the set $X$ is a complete poset. From the chain--continuousness of elements of \ffam\ it follows that they are isotone in $(X,\myle)$:
\beq\label{Fisot}
(x\myle y)\myimp(f(x)\myle f(y))\qquad\myll x,y\in X\ \myll f\in\ffam.
\eeq
By induction \fref{Fisot} implies the isotonicity of elements in \iter\ffam\ with respect to relation \myle:
\beq\label{ItFisot}
(x\myle y)\myimp(\varphi(x)\myle\varphi(y))\qquad\myll x,y\in X\ \myll\varphi\in\iter\ffam.
\eeq
In addition, as the family \ffam\ is commutative, so is the family \iter\ffam:
\beq\label{ItF-com}
\varphi\psi=\psi\varphi\qquad\myll\varphi,\psi\in\iter\ffam.
\eeq

Let us prove the equality
\beq\label{extF=extItF}
\ext{\iter\ffam}=\ext \ffam.
\eeq
Since $\ffam\subset\iter\ffam$, the embedding $\ext{\iter\ffam}\subset\ext\ffam$ takes place. Let us show the converse inclusion. Let $x'\in\ext\ffam$. We prove the inequality $x\myle\psi(x ')$ for all $\psi\in\iter\ffam$ by induction by the number of elements consisting $\psi$ (by the ``length'' of $\psi$). The induction base (for the ``length'' 1) follows directly from the choice of $x'$. Let the inequality $x\myle\varphi(x')$ hold for all $\varphi\in\iter\ffam$,  with the ``length'' not exceeding $k\in\NA$. Let $\psi'\in\iter\ffam$ have the ``length'' $k+1$. Therefore, $\psi'$ has the form $g\varphi'$, where $g\in\ffam$, $\varphi'\in\iter\ffam$ and the ``length'' of $\varphi'$ is equal to $k$. By the induction hypothesis, we have $x'\myle\varphi'(x')$. Then, by virtue of \fref{Fisot}, we obtain $g(x')\myle g\varphi'(x')=\psi'(x')$ and, in addition, by the choice of $x'$, we have $x'\myle g(x')$. By transitivity of \myle, from the last two relations we obtain the inequality $x'\myle\psi'(x')$. Since $\psi'$ is chosen arbitrarily, we have $x'\in\ext{\iter\ffam}$. Since the choice of $x'$ is arbitrary, the embedding $\ext\ffam\subset\ext{\iter\ffam}$ holds. This completes the proof of \fref{extF=extItF}.

Property \fref{extF=extItF} implies the fact that the set $\{\varphi(x): \varphi\in\iter\ffam\}$ is  directed for any $x\in\ext\ffam$. Indeed, using \fref{ItFisot}, \fref{ItF-com}, and \fref{extF=extItF}, we have
$$
\psi(x)\myle\psi(\varphi(x))=\varphi(\psi(x)),\qquad\varphi(x)\myle\varphi(\psi(x))\qquad\myll\psi,\varphi\in\iter\ffam\ \myll x\in\ext\ffam.
$$
Note that there are the same items in the right--hand side of these inequalities. Consequently, for all $\varphi,\psi\in\ffam$ the two-element sets $\{\varphi(x),\psi(x)\} $ has majorant. Using the transitivity of the relation $\myle$ we can extend this property by induction to an arbitrary finite subset of $X$. Thus, the set $\{\varphi(x): \varphi\in\iter\ffam\} $ is directed for any $x\in\ext\ffam$. Hence, since $(X,\myle)$ is complete poset, for any $x\in\ext\ffam$,
there exists $\sup\{\varphi(x):\varphi\in\iter{\ffam}\}$. The set $\ext\ffam$ is not empty, since $\bot\myle f(\bot)$ for all $f\in\ffam$. Thus, the right-hand side of \fref{fixfam} is well defined and non-empty set.

Let $x'\in\fix{\ffam}$. Then, obviously, $x'\myle f(x')$ and $ x'=\varphi(x')$ for all $f\in\ffam$ and $\varphi\in\iter\ffam$. Therefore, $x'\in\ext\ffam$ and $x'=\sup\{x'\}=\sup\{\varphi(x'):\varphi\in\iter{\ffam}\}$. That is, $x'\in\{\sup\{\varphi(x):\varphi\in\iter{\ffam}\}\mid x\in\ext\ffam\}$ and by virtue of the arbitrary choice of $x'$ we have
\beq\label{fixF-in-supit}
\fix\ffam\subset\{\sup\{\varphi(x):\varphi\in\iter\ffam\}\mid x\in\ext\ffam\}.
\eeq
Let us show the converse inclusion. Let $u$ be an arbitrary element of the right-hand side of \fref{fixF-in-supit}: $u=\sup\{\varphi(\bar x):\varphi\in\iter{\ffam}\}$, where $\bar x\in\ext\ffam$. Let $g\in\ffam$. By the condition of assertion (iii) $g$ is chain--continuous and, in addition, poset $X$ is chain--complete, and, so, strictly inductive (that is, every non-empty chain has the least upper bound). Therefore, for an arbitrary directed set $D\subset X$, holding in mind \cite[corollary 3]{Markow1976}, we have the equality
\beq\label{gsupst}
g(\sup D)=\sup\{g (x): x \in D\}.
\eeq
By virtue of the choice of $\bar x$, the sets $\bar D\mydef\{\varphi(\bar x): \varphi\in\iter\ffam\}$ and $\{g(x): x\in\bar D\} $ are directed; thus there exist the least upper bounds $\sup\bar D$ and $\sup\{g(x):x\in\bar D\}$, which, because of \fref{gsupst} and  the embedding $\{g(\varphi(\bar x)):\varphi\in\iter{\ffam}\}\subset\{\varphi(\bar x):\varphi\in\iter\ffam\}$, satisfy the relations
\beq\label{guleu}
g(u)=g(\sup\bar D)=\sup\{g(x): x\in\bar D\}\myle\sup\bar D = u.
\eeq

On the other hand, because of $\bar x$ choice, we have $\bar x\myle g(\bar x)$, and, by means of \fref{ItFisot} and \fref{ItF-com}, we get $\varphi(\bar x)\myle\varphi(g(\bar x))=g(\varphi(\bar x))$\ $\myll\varphi\in\iter{\ffam}$. In other words, $y\myle g(y)$\ $\myll y\in\bar D$. Then, once again using \fref{gsupst}, we get
\beq\label{gugeu}
u=\sup\bar D=\sup\{x:x\in\bar D\}\myle\sup\{g(x):x\in\bar D\}=g(\sup\bar D)=g(u).
\eeq
From \fref{guleu}, \fref{gugeu} it follows that $g(u)=u$. Due to the arbitrary choice of $g$ we have $u\in\fix\ffam$. Since $u$ is chosen arbitrarily we obtain the embedding
$$
\{\sup\{\varphi(x):\varphi\in\iter\ffam\}\mid x\in\ext\ffam\}\subset\fix{\ffam},
$$
which in conjunction with \fref{fixF-in-supit} gives the desired equality \fref{fixfam}.

We turn to the proof of \fref{botfixfam}. Let $x'\in\ext\ffam$. By definition of $\bot$ and \fref{ItFisot} it holds that $\varphi(\bot)\myle\varphi(x')$ $\myll\varphi\in\iter\ffam$. By the choice of $x'$ the set $D'\mydef\{\varphi(x'):\varphi\in\iter\ffam\}$ is directed in $(X,\myle)$. Since $(X,\myle)$ is a complete poset, there exists $\sup D'$. Of course, $\varphi(x')\myle\sup D'$ $\myll\varphi\in\iter\ffam$. From the last two inequalities we have
$$
\varphi(\bot)\myle\sup D'\qquad\myll\varphi\in\iter\ffam,
$$
that is, $\sup D'$ is a majorant of the set $\{\varphi(\bot):\varphi\in\iter\ffam\}$. By virtue of $\bot\in\ext\ffam$, the set $\{\varphi(\bot):\varphi\in\iter\ffam\} $ is ordered. Hence, there exists $\sup\{\varphi(\bot):\varphi\in\iter\ffam\}$ and, by definition of the least upper bound, we have the inequality
$$
\sup\{\varphi(\bot):\varphi\in\iter\ffam\}\myle\sup D'.
$$
Since $x'$ is chosen arbitrarily and due to \fref{fixfam}, we have
$$
\sup\{\varphi(\bot):\varphi\in\iter\ffam\}\in\fix\ffam.
$$
$$
\sup\{\varphi(\bot):\varphi\in\iter{\ffam}\}\myle u\qquad\myll u\in\fix\ffam.
$$
The relations imply the desired equality \fref{botfixfam}. Proof of theorem \ref{teo-limit-eq-fix-3} completed.

To substantiate corollary \ref{teo-limit-eq-fix-2}, we note that assertions (i) and (ii) of the corollary are a special cases of assertions (i) and (ii) of theorem \ref{teo-limit-eq-fix-3}, respectively. Equalities \fref{botfixf} and \fref{limit-eq-fix-2} are a special cases of equalities \fref{botfixfam} and \fref{fixfam}, respectively.


\begin{biblio}
\def\selectlanguageifdefined#1{
\expandafter\ifx\csname date#1\endcsname\relax
\else\language\csname l@#1\endcsname\fi}
\ifx\undefined\url\def\url#1{{\small #1}}\else\fi
\ifx\undefined\BibUrl\def\BibUrl#1{\url{#1}}\else\fi
\ifx\undefined\BibAnnote\long\def\BibAnnote#1{(#1)}\else\fi
\ifx\undefined\BibEmph\def\BibEmph#1{\emph{#1}}\else\fi

\bibitem{KraSub88e}
Krasovskii~N.~N., Subbotin~A.~I. \textit{Game-theoretical control problems}.
\newblock New York, Springer-Verlag Inc., 1988.

\bibitem{Chentsov76MS}
Chentsov~A.G. [On the game problem of pursuit at a given moment].{\it Matematicheskiy
sbornik.}
\newblock 1976.
\newblock V.~99(141), {\No}~3.
\newblock {pp.}~394--420. (in Russian)

\bibitem{Tarski1955}
Tarski~A. A lattice-theoretical fixpoint theorem and its applications. {\it
  Pacific Journal of Mathematics.}
\newblock 1955,
\newblock V.~5, {\No}~2,
\newblock {pp.}~285--309. 

\bibitem{Markow1976}
Markowsky~G. Chain--complete posets and directed sets with applications.
{\it  algebra universalis.}
\newblock 1976,
\newblock V.~6, {\No}~1,
\newblock {pp.}~53--68.

\bibitem{Ronse1994}
Ronse~Christian. Lattice--theoretical fixpoint theorems in morphological image
  filtering. {\it Journal of Mathematical Imaging and Vision.}
\newblock 1994,
\newblock V.~4, {\No}~1,
\newblock {pp.}~19--41.

\bibitem{Kantorovitch1939AM}
Kantorovitch~L. The method of successive approximation for functional equations. {\it Acta Mathematica.}
\newblock 1939,
\newblock V.~71, {\No}~1,
\newblock {pp.}~63--97.

\bibitem{BAREN1981E}
Barendregt~H.~P. {\it Lambda Calculus. Its Syntax and Semantics.}
\newblock Amsterdam New York Oxford, North-Holland Publishing Company, 1981.

\bibitem{Cohn1965}
Cohn~P.M. {\it Universal Algebra.}
\newblock New York, Harper and Row, 1965.
\end{biblio}

\newpage
\author{Серков Дмитрий Александрович, {\rm Институт математики и механики им. Н.Н. Красовского УрО РАН, 
Уральский федеральный университет им. Б.Н. Ельцина,
Россия, Екатеринбург, e-mail: serkov@imm.uran.ru}}
\title{Об аппроксимации совместных неподвижных точек}

\noindent
{\bf УДК 517.988.525, 517.977.8}

\ 

\maketitle{}

\begin{abstract} \begin{tabular}{p{0mm}p{139mm}}
&\noindent {\footnotesize \qquad 
Рассматривается семейство изотонных автоморфизмов частично упорядоченного множества. Известна теорема Альфреда Тарского о структуре множества совместых неподвижных точек таких автоморфизмов на полной решетке. Эта теорема имеет несколько обобщений (см., например, работы Марковского или Ронза), ослабляющих требования на прядковую структуру и модернизирующих соответствующим образом утверждение в части структуры множества совместных неподвижных точек. Вместе с тем, заметен недостаток утверждений типа теорем Канторовича или Клини, описывающих эти неподвижные точки как пределы последовательностей степней автоморфизмов. В статье приводятся условия на множество и семейство его автоморфизмов, при которых  итеративные последовательности элементов рассматриваемого семейства аппроксимируют множество совместных неподвижных точек. Данный результат развивает в конструктивном направлении упомянутые утверждения Тарского, Марковского и Ронза.

\qquad\keywords{совместные неподвижные точки, итерационный предел}}
\end{tabular}\end{abstract}

\begin{flushright}
{ \it Received ?? ?? ??}
\end{flushright}


\begin{thebibliography}{1}
\def\selectlanguageifdefined#1{
\expandafter\ifx\csname date#1\endcsname\relax
\else\language\csname l@#1\endcsname\fi}
\ifx\undefined\url\def\url#1{{\small #1}}\else\fi
\ifx\undefined\BibUrl\def\BibUrl#1{\url{#1}}\else\fi
\ifx\undefined\BibAnnote\long\def\BibAnnote#1{(#1)}\else\fi
\ifx\undefined\BibEmph\def\BibEmph#1{\emph{#1}}\else\fi

\bibitem{KraSub88e}
Krasovskii~N.~N., Subbotin~A.~I. Game-theoretical control problems.
\newblock Springer-Verlag New York, Inc., 1988.
\newblock {с.}~517.

\bibitem{Chentsov76MS}
Ченцов~А.Г. Об игровой задаче сближения в
  заданный момент времени~// Математический
  сборник.
\newblock 1976.
\newblock Т. 99(141), {\No}~3.
\newblock {С.}~394--420.

\bibitem{Tarski1955}
Tarski~A. A lattice-theoretical fixpoint theorem and its applications~//
  Pacific Journal of Mathematics.
\newblock 1955.
\newblock Т.~5, {\No}~2.
\newblock {С.}~285--309. URL:  http://projecteuclid.org/euclid.pjm/1103044538.

\bibitem{Markow1976}
Markowsky~George. Chain--complete posets and directed sets with applications~//
  algebra universalis.
\newblock 1976.
\newblock Т.~6, {\No}~1.
\newblock {С.}~53--68.

\bibitem{Ronse1994}
Ronse~Christian. Lattice--theoretical fixpoint theorems in morphological image
  filtering~// Journal of Mathematical Imaging and Vision.
\newblock 1994.
\newblock Т.~4, {\No}~1.
\newblock {С.}~19--41.

\bibitem{Kantorovitch1939AM}
Kantorovitch~L. The method of successive approximation for functional
  equations~// Acta Mathematica.
\newblock 1939.  December.
\newblock Т.~71, {\No}~1.
\newblock {С.}~63--97.

\bibitem{BAREN1981E}
Barendregt~H.~P. Lambda Calculus. Its Syntax and Semantics.
\newblock Amsterdam New York Oxford: North-Holland Publishing Company, 1981.
\newblock Т.~103 {из} \BibEmph{Studies in Logic and Foundations of
  Mathematics}.

\bibitem{Cohn1965}
Cohn~P.M. Universal Algebra.
\newblock New York: Harper and Row, 1965.
\end{thebibliography}
\end{document}